\documentclass[11pt]{article}

\usepackage[dvips]{epsfig}
\usepackage{graphicx}
\usepackage{amssymb,graphics,amsmath,amsthm,amsopn,amstext,amsfonts}
\usepackage{color}
\newcommand{\supp}{ \mathrm{supp} }

\newcommand{\argmax}{\operatornamewithlimits{\arg\max}}
\newcommand{\argmin}{\operatornamewithlimits{\arg\min}}

\newcommand{\Jcal}{\mathcal J}

\newcommand{\mycut}[1]{{}}
\newtheorem{theorem}{Theorem}[section]

\newtheorem{definition}{Definition}[section]

\usepackage{verbatim}
\usepackage{listings}
\usepackage{graphicx}
\usepackage{subcaption}

\providecommand{\keywords}[1]{\textbf{\textit{Keywords.}} #1}

\date{}

\begin{document}

\title{\bf{A Survey on Compressive Sensing: Classical Results and Recent Advancements}}

\author{Ahmad Mousavi\footnote{Institute for Mathematics and its Applications, University of Minnesota, Minneapolis, MN, USA }\ ,  Mehdi Rezaee\footnotemark[2]
  \ and   Ramin Ayanzadeh\footnote{Department of Computer Science and Electrical Engineering, University of Maryland Baltimore County, Baltimore, MD 21250, USA
   \newline {Emails: amousavi@umn.edu, and \{rezaee1,  ayanzadeh\}@umbc.edu.} 
  } 
}
\maketitle

\begin{abstract}
Recovering sparse signals from linear measurements has demonstrated outstanding utility in a vast variety of real-world applications.  Compressive sensing is the topic that studies the associated raised questions for the possibility of a successful recovery. 
This topic is well-nourished and numerous results are available in the literature. However, their dispersity makes it challenging and time-consuming for readers and practitioners to quickly grasp its main ideas and classical algorithms, and further touch upon the recent advancements in this surging field. Besides, the sparsity notion has already demonstrated its effectiveness in many contemporary fields. Thus, these results are useful and inspiring for further investigation of related questions in these emerging fields from new perspectives. In this survey, we gather and overview vital classical tools and algorithms in compressive sensing and describe significant recent advancements. We conclude this survey by a numerical comparison of the performance of described approaches on an interesting application. 
\end{abstract}
\keywords compressive sensing, $\ell_p$ recovery, greedy algorithms

\section{Introduction}
In traditional sensing, i.e., uniform sampling, we need to densify measurements to obtain a higher-resolution representation of physical systems. But in  applications like multiband signals with wide spectral ranges, the required sampling rate may exceed the specifications of the best available analog-to-digital converters (ADCs) \cite{rani2018systematic}. On the other hand, measurements obtained from linear sampling methods  approximately carry a close amount of information, which makes them reasonably robust; for example,  for packages lost in data streaming applications \cite{pudlewski2012compressed}.
Compressive sensing linearly samples sparse signals at a rate  lower than the traditional rate provided by the Nyquist-Shannon sampling theorem \cite{yaroslasky2015compressed}. 

Recovering a sparse signal from a  linear measurement is of high interest to preserve more storage, to have  less computation, energy, and communication time, and to propose effective  data compression and transmission methods \cite{ayanzadeh2019compressive-geo, bai2020high, baraniuk2007compressed,  boche2015compressed, candes2006compressive, EldarKutyniok_book2012, FoucartRauhut_book2013, kutyniok2013theory, rani2018systematic, RishGrab_book2014, sankaranarayanan2020compressive, wang2020lightamc, zhang2015survey, zhao2018sparse}. 
There are also specialised motivations behind this interest; for example, the recovery process in compressive sensing requires prior knowledge about original signals, which adds reasonable levels of security and privacy to the compressive sensing based data acquisition framework \cite{zhang2018secure}.
 Compressive sensing-based face recognition techniques are invariant to rotation, re-scaling and data translation \cite{rani2018systematic}. Besides, compressive sensing has demonstrated outstanding numerical stability in terms of not only noisy measurements but also quantization errors \cite{saab2018quantization}. In machine learning, compressive sensing has demonstrated to boost the performance of pattern recognition methods.

The sparsity assumption as the basis of compressive sensing is not practically   restrictive because it is empirically observed that desired signals in an extensive range of applications are sparse, possibly after a change of basis though. However, blind compressive sensing techniques are universal in terms of transfer domains, namely, there is no need to know  the sparsifying basis while taking linear measurements \cite{gleichman2011blind}.  Sparsity notion directly leads to NP-hard problems and thus the original problem of compressive sensing is computationally intractable; see Section \ref{sec:original} for the definition of this problem. Two main strategies to tackle this obstacle are  convex or nonconvex relaxations and greedy algorithms; see Sections \ref{sec:lp_recovery} and \ref{sec:greedy}. 
There are several other effective approaches that have different ideas than these two categories. Hence, this survey  attempts to briefly characterize these promising studies as well. 

To relax an $\ell_0$-quasi-norm based sparse optimization problem, $\ell_0$-quasi-norm is often replaced by $\ell_p$-(quasi)-norm  with $p>0$. Then, it is investigated that under which conditions on problem parameters and for which $p$'s both original and approximate problems uniquely obtain the same solution.  The efficiency of $p=1$ is well-documented \cite{CarWangX_TSP10,  Cai_2009, Candes_TIT2005, Donoho_SR2001, DonohoElad_PNAS2003, DonohoTanner_PNAS2005, Foucart_ACHA10}. But since this choice does not lead into strictly convex programs, the solution uniqueness shortcoming needs to be addressed \cite{mousavi_2017, Zhang_Yin_Cheng}.   This case converts to a linear program and since those matrices involved in sparse optimization often inherit large dimensions, well-known methods  including simplex and interior-point methods or specified algorithms are applied \cite{burger2013adaptive,salahi2007mehrotra}. The shape of unit balls associated with $0 < p < 1$ encourages examining this case as well. It turns out that  the  recovery results for certain values of $0 < p < 1$ are more robust and stable compared to  $p=1$ \cite{Chartrand_2007, Chartrand_2009}. Despite these favorable  theoretical results, nonconvex relaxations ask for a global minimizer, which is an intractable task. One numerical way to bypass this is by utilizing classical schemes that obtain a local minimizer with an initial point sufficiently close to a global minimizer. For example, the solution of least squared is an empirical suggestion for this initialization \cite{Chartrand_2007}, although there is no guarantee that such a solution is close enough. 
Nonetheless, a recent  theoretical study proves that the main representative  problems utilized in the realm of sparse optimization, such as the generalized basis pursuit and LASSO, almost always return a full-support solution for $p>1$ \cite{Shen_2018}.

Greedy algorithms in our second class have  low computational complexity, especially for relatively small sparsity levels, and yet they are effective. A well-received algorithm in this class which plays a key role is the orthogonal matching pursuit (OMP).  To directly solve the original problem of compressive sensing, it exploits the best local direction in each step and adds its corresponding index to the current support set. Then, it estimates the new iteration as the orthogonal projection of the measurement vector onto the subspace generated by columns in the current support set \cite{Tropp_TIT2004, Tropp_TIT2007}. There is an extensive literature on the capability of this algorithm in identifying the exact support set that vary based on the required number of steps and accounting for noise in the measurement vector. For example, a sufficient condition for recovering an $s$ sparse vector from an exact linear measurement asks for an $(s+1)$st restricted isometry constant strictly smaller than $(\sqrt{s}+1)^{-1}$  \cite{Mo_2012}. This upper bound improves to the necessary and sufficient bound of $1/\sqrt{s+1}$.  A result with the same spirit for the noisy measurement is available that imposes an extra assumption on the magnitude of the smallest nonzero entry of a desired sparse vector \cite{Wen_2013}. Further, stability is also achievable with the presence of noisy measurements if the number of required steps is increased  \cite{Zhang_2011}. Nevertheless, there are other algorithms in this class that have the same spirit but attempt to resolve issues of the OMP, such as stagewise OMP \cite{DonohoTDS_TIT12} and multipath matching pursuit \cite{KwonWS_TIT14}. %
For example, many algorithms allow more indices to enter the current support set  for picking a correct index in each step such as the compressive sampling matching pursuit \cite{needell2008iterative} and  subspace pursuit \cite{dai2009subpsace}. Other effective greedy algorithms extending the OMP are simultaneous, generalized, and grouped  OMP \cite{SwirszczAL_NIPS09,TroppGS_SP06, WangKS_TSP2012,wen2017novel}. Nonnegative and more general constrained versions of the OMP obtain similar recovery results based on  variations of  the restricted isometry and the orthogonality constants \cite{Bruckstein&Elad&Zibulevsky_TIT2008,lin2014stable,shen_2019}. These algorithms are also generalized for block sparse signals \cite{eldar2010block,wen2018optimal}.

The rest of paper is organized as follows. In Section \ref{sec:original}, we first introduce the original problem of compressive sensing and its extensions and  variations.  Section \ref{sec:lp_recovery} is devoted to $\ell_p$ recovery which covers using convex and nonconvex relaxations. In Section \ref{sec:greedy}, we study several well-received greedy algorithms. Finally, we report numerical performance results in Section \ref{sec:numerical_results}.

\textit{Notation}.
For a natural number $N$, we let $[N]:=\{1,2,\dots,N\}$. For a set $S\subseteq [N]$, its number of elements  and complement are denoted  by $\text{card}(S)$ (or $|S|$) and $S^c$, respectively. 

For a vector $x\in \mathbb R^N$ and $p>0$, let $\|x\|_p:=(\sum_{i\in [N]} |x_i|^p)^{1/p}$. For $S=\{i_1,i_2\dots, i_{|S|}\} \subseteq [N]$, we let $x_S:=z\in \mathbb R^N$ such that $z_i=x_i$ for $i\in S$ and $z_i=0$ whenever $i\in S^c$, that is, $x_S$ denotes the coordinate projection of $x$ onto the subspace generated by $\{e_{i_1}, e_{i_2}, \dots, e_{i_{|S|}}\}$. Given $x\in \mathbb R^N$,  the hard thresholding  operator of order $s$ is defined as $H_s(x):=x_S$ where $S$ contains  $s$ largest absolute entries of $x$ (notation $L_s(x)$  denotes this index set). Further, let $x_+:=\max (x,0)$. A vector is called $s$-sparse if it has at most $s$ nonzero entries, i.e., its sparsity level is bounded by $s$. We use $\langle ., . \rangle$ to denote the standard  inner product between two vectors. We employ $\textbf{1}_N$ to denote a vector of size $N$ with each entry equal to $1$, and simply let $0$  denote a zero vector, which its dimension should be clear based on the discussed content. 
For a matrix $A\in \mathbb{R}^{m\times N}$, we denote its  transpose by $A^T$, its adjacent by $ A^*$ and its Frobenius norm by $\|A\|_F$. We denote the $j$th column of $A$ by $A_j$ and, given a subset $S\subseteq [N]$, a corresponding column submatrix is defined as  $A_S:=[A_{j_1}, A_{j_2}, \dots , A_{j_{|S|}}]$. The null and range spaces of this matrix are denoted by $\mathcal Ker(A)$ and $\mathcal R(A)$, respectively. The symbol $\mathcal O$ stands for the standard big O notation.

\section{Compressive Sensing Problem and Its Extensions and  Variations} \label{sec:original}

Compressive sensing (CS), also known as compressed sensing or sparse sampling, is established based on sparsity assumption. Hence,  we start with preliminary definitions of sparsity and, a closely related concept, compressibility.  
\begin{definition}
For a vector $x\in \mathbb R^N$, let $\supp(x):=\{i\in [N] \,  | \, x_i\ne 0\}$ and $\|x\|_0:=\mathrm{card}(\supp(x))$. In particular, for $s\in [N]$, this vector is called $s$-sparse if $\|x\|_0\le s$.
\end{definition}
\begin{definition}
For a vector $x\in \mathbb R^N$,  $s\in [N]$, and $p>0$, let
\[
{\sigma_s(x)}_p:=\min_{z\in \mathbb R^N} \|x-z\|_p  \quad \mathrm{subject ~ to} \quad \|z\|_0\le s.\]
This vector is called $s$-compressible if ${\sigma_s(x)}_p$ is small for some $p>0$. This vector is sometimes called nearly $s$-sparse in $\ell_p$-norm.
\end{definition}
It is easy to show that ${\sigma_s(x)}_p=\|x-x_S\|_p$, where $S\subseteq [N]$ contains all the  $s$ largest absolute entries of $x$ (consequently, ${\sigma_s(x)}_p=0$ when $x$ is $s$-sparse). Roughly speaking, a vector is called sparse if most of its entries are zero and it is called compressible if it is well-approximated by a sparse vector. 
Sparsity as a prior assumption for a desired vector is  consistent with diverse applications  because it is often doable to employ a change of basis technique for finding sparse representations in a transform domain.  For instance, wavelet, Radon, discrete cosine, and Fourier transforms are well-known as suitable and efficient choices for natural, medical, and digital images, and speech signals, respectively.

The original problem of compressive sensing aims to recover a sparse signal $ {x} $ from a linear measurement vector $ {y}=  A {x}$, where $  A \in \mathbb{R}^{m\times N}$  (with $m \ll N$) is the so called measurement, coding or design matrix. So, the ultimate goal of compressive sensing  is reasonably formulated as follows: 
\begin{equation}   \tag{CS}  \label{pr: initial_cs_problem}
\min_{  x \in \mathbb R^N} \
\| {x}\|_0
\quad \text{subject to} 
\quad     A {x}= {y}.
\end{equation}
To incorporate noisy measurements, one can study this problem: 
\begin{equation*}   \label{pr:noisy_ initial_cs_problem}
\min_{  x \in \mathbb R^N} \
\| {x}\|_0
\quad \text{subject to} 
\quad    \|Ax-y\|_2\le \eta,
\end{equation*}
where $\eta>0$ controls the noise level.
These problems have demonstrated to revolutionize several real-world applications in both science and engineering disciplines, including but not limited to signal processing, imaging, video processing, remote sensing, communication systems, electronics, machine learning, data fusion, manifold processing, natural language and speech processing, and processing  biological signals \cite{marques2018review,qiasar2013compressive,qin_2018_sparse,rontani2016compressive,zhang_2019_electrocardiogram}. 

In a broad range of real-world applications, we utilize ADCs to map the real-valued measurements of physical phenomena (over a potentially infinite range) to discrete values (over a finite range). As an illustration, when we use $b$ bits for representing the measurements digitally, the quantization module in an ADC maps measurements to one of the $2^{b}$ distinct values that introduces error in measurements. Quantization module plays a bottleneck role in restricting the sampling rate of the ADCs because  the maximum sampling rate decreases exponentially when the number of bits per measurement increases linearly. Quantization module also is the main source of energy consumption in ADCs. Single-bit (or 1-bit) compressive sensing enables us to reduce the number of bits per measurements to one and introduces a proper model for successful recovery of the original signal. This extreme quantization approach only retains the sign of  measurements, i.e., we have $y_i\in \{1,-1\}$ for all $i\in [m]$, which results in a significantly efficient, simple, and fast quantization \cite{boufounos20081, jacques2013robust}. 

A recent generalization of  (\ref{pr: initial_cs_problem}) is the so called compressive sensing with matrix uncertainty.  This uncertainty finds two formulations that incorporate measurement errors as well. The first formulation is the following:
\begin{equation*}    
\min_{  x \in \mathbb R^N, \ E \in \mathbb R ^{m\times N}} \
\|(A+E)x-y\|^2_2+\lambda_E \|E\|^2_F+\lambda \| {x}\|_0,
\end{equation*}
where the matrix $E$ is the perturbation matrix. The second one is  the following:
\begin{equation*}    
\min_{  x \in \mathbb R^N,\ d \in \mathbb R ^{r}} \
\bigg \|(A^{(0)}+\sum_{i\in [r]}d_iA^{(i)})x-y\bigg \|^2_2+\lambda_d \|d\|^2_2+\lambda \| {x}\|_0,
\end{equation*}
where the measurement matrices $A^{(i)}$'s for $i=0,1,\dots, r$ are known and the unknown vector $d$ is the uncertainty vector. 
In quantized compressed sensing, we have either $y=\mathcal QAx$ or $y=\mathcal Q (Ax+v)$ where $v$ is the noise such that $v\sim \mathcal N(0, \sigma^2I)$ and $\mathcal Q:\mathbb R^m \rightarrow \mathcal A \subseteq \mathbb R^m$ is the set-valued quantization function. Well-studied examples of such function, map $Ax$ or $Ax+v$ into $\mathcal A=\{x\in \mathbb R^m \ | \ l\le x \le u\}$ or $\mathcal A=\{+1,-1\}^m$ \cite{jacques2013robust,xu2018quantized,zymnis2009compressed}.
The former case leads to:
\begin{equation*} 
\min_{  x \in \mathbb R^N} \
\| {x}\|_0
\quad \text{subject to} 
\quad     l\le Ax+v\le u,
\end{equation*}
and the latter one yields the following problem:
\begin{equation*} 
\min_{  x \in \mathbb R^N} \
\| {x}\|_0
\quad \text{subject to} 
\quad     y=\mbox{sign}(Ax+v).
\end{equation*}
In some recent applications, we may have nonlinear measurements that encourage the so called cardinality constrained optimization:
\begin{equation*} \label{eqn:general_cardinality}
\underset{x \in \mathbb R^N}{\text{min}} \ f(x) \ \quad \text{subject to} \quad \|x\|_0\le s.
\end{equation*}
Two important examples are $f(x)=\|Ax-y\|^2$ (motivated by linear measurements $y=Ax$) and $f(x)=\sum_i^m (y_i-x^TA^{(i)}x)^2$ (proposed for quadratic measurements of  $y_i=x^TA^{(i)}x$ for $i=1,2,\dots, m$ and symmetric matrices $A^{(i)}$'s) \cite{shechtman2011sparsity, szameit2012sparsity}.
Many algorithms have been proposed for solving this cardinality constrained problem with specific and general smooth functions \cite{ beck2013sparse,BeckE_SIOPT13,Beck_2015,BKanzowS_SIOP16,shechtman2014gespar}.
In matrix setting, this problem converts to rank minimization, which finds numerous applications \cite{FazelHBody_ACC01}. Although we listed some vital variations of compressive sensing above, our focus in this survey is on the original problem to motivate the important tools, algorithms and results in compressive sensing for a linear measurement possibly contaminated with some error or noise.

\section{$\ell_p$ Recovery with $0<p\le 1$: Main Formulations and Results} \label{sec:lp_recovery}
Since the  problem (\ref{pr: initial_cs_problem})  is  NP-hard \cite{Muthukrishnan_2005}, it is not tractable and it must be handled indirectly. Because $\ell_p$-norm approximates $\ell_0$-norm as $p$ goes to zero, one approach to tackle (\ref{pr: initial_cs_problem}) is  exploiting the so called Basis Pursuit (BS) problem:
\begin{equation}  \tag{BP} \label{pr: Basis_Pursuit}
 \min_{  x \in \mathbb R^N} \
\| {x}\|_1
\quad \text{subject to} 
\quad    A {x}= {y}.
\end{equation}
Then, one investigates  under which conditions (\ref{pr: initial_cs_problem}) and (\ref{pr: Basis_Pursuit}) appoint the same solution.
This equivalence property occurs when  the solution  is highly sparse and the measurement matrix has sufficiently small  mutual incoherence (defined below) \cite{DonohoElad_PNAS2003}. 
\begin{definition}
For a matrix $A\in \mathbb{R}^{m\times N}$, the mutual coherence is defined as 
\[
\mu(A):=\max_{i\ne j} \frac{|\langle A_i, A_j \rangle|}{\|A_i\|_2.\|A_j\|_2},
\]
where $A_i$ denotes the $i$th column vector of the matrix $A$.
\end{definition}
The mutual coherence  simply seeks the largest correlation between  two different columns. The following global 2-coherence is a generalization of the mutual coherence ($\mu = \nu_1$):
\begin{definition}
For a matrix $A\in \mathbb{R}^{m\times N}$, its $k$th global 2-coherence is defined as below:
\[
\nu_k (A):=\max_{i\in [N]}
\ \max_{\Lambda \subseteq [N]\setminus \{i\}, |\Lambda|\le k} 
\Bigg(\sum_{j\in \Lambda }\frac{\langle A_i, A_j \rangle^{2}}{\|A_i\|^2_2.\|A_j\|^2_2} \Bigg)^{1/2}.
\]
\end{definition}
This tool is useful for a successful recovery of a weak version of orthogonal matching pursuit \cite{yang2014new}, while we focus on the standard mutual coherence here to convey the main idea.

\begin{theorem}\cite{Elad_2002} \label{thr:coh}
Assume that $\|x\|_0<0.5(1+\mu(A)^{-1})$. Then, problems \ref{pr: initial_cs_problem}  and \ref{pr: Basis_Pursuit} uniquely obtain the signal $x$.
\end{theorem}
Sparsity is an inevitable assumption in the recovery process. This assumption accompanied by  other conditions on a measurement  matrix such as  the mutual coherence or restricted isometry property (RIP) guarantee a successful recovery through the   problem \ref{pr: Basis_Pursuit}. For example, the bound in Theorem \ref{thr:coh} is improved using the RIP tool defined below.
\begin{definition}
Let $A$ be an $m\times N$ matrix. Then, this matrix has the restricted isometry property of order $s$ provided that there exists $\delta_s\in (0,1)$  such that 
\[
(1-\delta_s)\|x\|_2^2\ \le \ \|Ax\|_2^2 \ \le \ (1+\delta_s)\|x\|_2^2 \quad \quad \quad \forall x; \ \|x\|_0\le s. 
\]
\end{definition}
This definition demands each column submatrix $A_S$ with $\mbox{card}(S)\le s$ to have singular values in $[1-\delta_s, 1+\delta_s].$ Consequently, such  submatrices have linearly independent columns if $\delta_s<1$. 
Since this definition involves all the $s$-tuples of columns, it is more rigorous  than the mutual coherence ($\delta_2=\mu$). Thus, the RIP tool derives to better upper bounds on the sparsity level of a vector to be recovered, that is, less sparse vectors are properly handled. 
\begin{theorem}\cite{Mo_2011}
Problems \ref{pr: initial_cs_problem} and \ref{pr: Basis_Pursuit} uniquely obtain the same $s$-sparse solution if $\delta_{2s}<0.4931$.
\end{theorem}

There are numerous similar results for recovering  $s$-sparse  signals in  compressive sensing literature based on the RIP constants, which  are often in the form of $ \delta_{ks}\le \delta$ for some numbers $k>0$ and $\delta\in (0,1)$. For example, the \ref{pr: Basis_Pursuit} recovers all the $s$-sparse signals if $\delta_s<1/3$ and this bound is  sharp \cite{Cai_2013}. Since the condition $\delta_{ks}\le \delta$ is not practically verifiable due to its computational complexity \cite{Tillmann_2014}, a breakthrough in the compressive sensing field took place when  random matrices, such as Gaussian and Bernoulli, demonstrated their capability to satisfy such appealing inequalities with a high probability. 
In fact, for an $m\times N$ random matrix, where each entry is independently drawn from Gaussian or Bernoulli distributions, we have  $\delta_{s}\le \delta$ for  $m \ge C \delta^{-2} s \ln(\frac{eN}{s})$ in which $C >0$ does not depend on $s, m$ and $N$, and $e$ is the natural number \cite{Baraniuk_2008}. It must be mentioned that most results for constructing deterministic matrices with the RIP property depend on the mutual coherence, leading to $m\ge cs^2$ with $c>0$. This bound is not practical because the number  of measurements $m$ scales  quadratically in sparsity level $s$ for deterministic matrices versus linearly for random matrices. This confirms the desire of using random matrices in compressive sensing. On the other hand, there are recent promising results toward the practicality of deterministic matrices in terms of not only the compression ratio of $m/N$ but also the complexity of the recovery process 
\cite{amini2011deterministic, candes2007sparsity, gu2019deterministic, nguyen2013deterministic}. For example, deterministic sparse coding matrices have demonstrated to outperform both binary tree recovery and binary $\ell_1$-norm recovery methods in the so called binary compressive sensing \cite{nakarmi2012bcs}.
The majority of analyses within  compressive sensing relies on the RIP but this tool has minor drawbacks that are partially rectified via the NSP  defined below. 
\begin{definition}
Let $A$ be an $m\times N$ matrix. This matrix has the null space property of order $s$ if  
\[
\|v_{S}\|_1<\|v_{S^c}\|_1; \quad  \forall v\in \mathcal{K}er(A)\setminus\{0\} \ \mathrm { and } \ \forall S\subseteq [N]\ \mathrm{ with } \ |S|\le s.
\] 
\end{definition}
In principle, the measurement  pair $(A,y)$ carries all the required information for the recovery process. Since linear systems  $Ax=y$ and $PAx=Py$ have identical solution sets for a nonsingular matrix $P$, the measurement pair $(PA,Py)$ has the same information but it  is numerically a better choice for a suitable conditioning matrix $P$.  Although the RIP constants of matrices $A$ and $PA$ can vastly differ \cite{Y_Zhang_2013}, the NSP is preserved if either one holds it. In addition, this property is  necessary and sufficient  for the  solution uniqueness  of problem \ref{pr: Basis_Pursuit} (in the uniform sense). 
\begin{theorem}\cite{FoucartRauhut_book2013}
Every $s$-sparse vector $x$ is the unique solution to \ref{pr: Basis_Pursuit} with $y=Ax$ if and only the measurement matrix $A$ holds the NSP of order $s$. 
\end{theorem}
The above theorem implies that under the NSP of order $s$, the convex program \ref{pr: Basis_Pursuit} as a matter of fact solves the NP-hard problem \ref{pr: initial_cs_problem}.

The Range Space Property (RSP) is another  useful condition in compressive sensing given by Zhao \cite{Zhao_2013}.
\begin{definition}
The matrix $A^T$ has the range space property of order $s$ if for any disjoint subsets $S_1$ and $S_2\subseteq [N]$ with $|S_1|+|S_2|\le s$ there exists a vector $\eta\in \mathcal{R}(A^T)$ such that 
\[
\eta_i=1 \quad \forall i\in S_1; \quad \eta_i=-1 \quad \forall i\in S_2 \quad \mbox{ and }\quad  \|\eta_{{(S_1\cup S_2)}^c}\|_{\infty}<1.
\]
\end{definition}
\begin{theorem}\cite{Zhao_2013}
Every $s$-sparse vector $x$ is exactly recovered by the problem \ref{pr: Basis_Pursuit} with $y=Ax$ if and only if $A^T$ has the RSP of order s.
\end{theorem}
The matrix $A^T$ has the range space property of order $s$ under several assumptions such as $s<0.5(1+\mu(A)^{-1}), \ \delta_{2s}(A)<\sqrt{2}-1$ and the NSP of order $2s$ \cite{Zhao_2013}. Further, an extended version of this property is useful to provide a similar result to the above theorem for the nonnegative sparse recovery. 

So far, we have only discussed noiseless measurements in this section, although a practical situation definitely imposes investigating noisy scenario as well. Taking into account noisy measurements yields in the following  Quadratically Constrained Basis Pursuit (QCBP) problem (after replacing $\ell_0$ by $\ell_1$):
\begin{equation*}  \tag{QCBP} \label{pr: Quadratically_Constrained_Basis_Pursuit}
 \min_{  x \in \mathbb R^N} \
\|{x}\|_1
\quad \text{subject to} 
\quad  \|Ax-y\|_2\le \eta.
\end{equation*}
\begin{theorem}\cite{FoucartRauhut_book2013}
Given $x\in \mathbb R^N$ and a  matrix $A\in \mathbb{R}^{m \times N}$ such that $\delta_{2s}<4/\sqrt{41}\approx 0.6246$, every minimizer $x^*$ of \ref{pr: Quadratically_Constrained_Basis_Pursuit} satisfies 
$\|x^*-x\|_2\le Cs^{-1/2}{\sigma_s(x)}_1 +D\eta$ with $C>0$ and $D>0$.
\end{theorem}
The following so called Least Absolute Shrinkage and Selection Operator (LASSO) is practically more efficient \cite{Wright_Nowak_Figueiredo,Wainwright_2009,WenYinGoldfarbZhang_SJSC2010}: 
\begin{equation*}  \tag{LASSO} \label{pr: LASSO}
 \min_{  x \in \mathbb R^N }\
\quad  \frac{1}{2}\|  y-   A {x}\|^2_2+\lambda \| {x}\|_1.
\end{equation*}
There is a trade-off between feasibility and sparsity in the nature of this problem, which is controlled by the regularization parameter $\lambda$. Not only the \ref{pr: LASSO}  has been extensively studied \cite{Tibshirani_RSS1996}, but also there are various results for it in terms of its effectiveness in sparse optimization. For instance, we bring the following theorem for deterministic measurement matrices. 

\begin{theorem}\cite{Wainwright_2009}
Let $y=Ax+\textbf{e}$ where $x$ is supported on $S\subseteq [N]$ with $|S|\le s$ and $\textbf{e}$ is a zero-mean additive observation noise.
Assume that the measurement matrix $A$ satisfies $\|A^T_{S^c}A_S(A^T_S A_S)^{-1}\|_{\infty}<1$, $\lambda_{\min} (A^T_S A_S)>m$ and $\max_{j\in S^C} \|A_j\|_2\le \sqrt{m}$. Further, assume that $N=\mathcal{O}(\exp(m^{\gamma})), \ s=\mathcal{O}(m^\alpha)$, $\min_{i\in S} x_i>1/m^{\frac{1-\beta}{2}}$ with $0<\alpha+\gamma< \beta<1$. For $\lambda=m^{\frac{1-\delta}{2}}$ such that $\delta \in (\gamma, \beta-\alpha)$, the \ref{pr: LASSO} problem recovers the sparsity pattern $S$ with probability $1-\exp(-cm^\delta)$ for a constant $c>0$.
\end{theorem}
For the case of random measurement matrices, this problem  needs a sample size $m > 2s \ln(N-s)$ to achieve exact recovery with a high probability. This probability  converges to one for larger problems \cite{Wainwright_2009}. An efficient  algorithm for solving the  BP utilizes a sequence of the LASSO problems  \cite{yin2008bregman}.
The following so called Basis Pursuit Denoising (BPD) minimizes the feasibility violation while implicitly bounding the sparsity level with a parameter $\tau \ge 0:$
\begin{equation*}  \tag{BPD} \label{pr: Basis Pursuit Denoising}
 \min_{  x \in \mathbb R^N} \ 
 \|  y-   A {x}\|_2
\quad \text{subject to} 
\quad  \| {x}\|_1 \le \tau.
\end{equation*}
There are known relations among the optimums of  the \ref{pr: Quadratically_Constrained_Basis_Pursuit}, \ref{pr: LASSO}, and \ref{pr: Basis Pursuit Denoising} problems, which can be found in \cite{FoucartRauhut_book2013}.
The Dantzig Selector (DS) problem arises in several statistical applications, so it has been also employed in sparse recovery \cite{Cai_2009,Candes_2007}:
\begin{equation*}  \tag{DS} \label{pr: Dantzig Selector}
 \min_{  x \in \mathbb R^N} \ 
 \|  x \|_1
\quad \text{subject to} 
\quad  \| A^T(Ax-y)\|_\infty \le \sigma. 
\end{equation*}
This  problem    manages noisy measurements and  reduces to  a linear programming problem. To bring a pertaining result next, we need the following restricted  orthogonality quantity.
\begin{definition}
Given a matrix $A\in \mathbb{R}^{m\times N}$, its $s,s'$-restricted orthogonality (RO) constant $\theta_{s,s'}$ is defined as the smallest $\theta>0$ such that 
\[
|\langle Ax, Ax' \rangle|\le \theta \|x\|_2 \|x'\|_2 \quad \forall x~\mathrm{ and }~x'; \ \|x\|_0\le s~\mathrm{and} \ \|x'\|_0\le s'.
\]
\end{definition}
 According to \cite{Cai_2009}, for $x\in \mathbb R^N,$ and a noisy $y$ such that $\|A^T(Ax-y)\|_\infty \le \sigma$ and $\delta_{1.5s}+\theta_{s,1.5s}<1$, an optimal solution $x^*$ to \ref{pr: Dantzig Selector} obeys 
$
  \|x^*-x\|_2\le Cs^{\frac{1}{2}}\sigma+Ds^{\frac{-1}{2}}{\sigma_s(x)}_1,$ 
  where $C$ and $D$ are two constants depending on $\delta_{{1.5}s}$ and $\theta_{s,1.5s}$.
  In particular, if $x$ is an $s$-sparse vector, then $\|x^*-x\|_2\le Cs^{\frac{1}{2}}\sigma.$
A very recent result related to the recovery of  this problem as reported in Corollary 5.6.1 of \cite{zhao2018sparse} is brought below.
\begin{theorem}\cite{zhao2018sparse}
Let $A\in \mathbb R^{m\times N}$ with $m<N$ and $\text{rank}(A)=m.$  Consider the \ref{pr: Dantzig Selector} problem.
If $A^T$ satisfies the RSP of order $k$, then for any $x\in \mathbb R^n$, there is a solution $x^*$ of the \ref{pr: Dantzig Selector} approximating $x$ with error
\[
\|x-x^*\|_2\le \gamma\{2\sigma_k(x)_1+\left(
\|A^T(Ax-y)\|_\infty-\sigma
\right)_+
+c_A\sigma+c_A\|A^T(Ax-y)\|_\infty\},
\]
where $\gamma$ is a constant and $c_A$ is a constant given as 
\[
c_A= \max_{G\subseteq \{1,\dots,n\}, |G|=m} \|A^{-1}_G(AA^T)^{-1} A\|_{\infty \rightarrow 1},
\]
where $A_G$ stands for all positive $m\times m$ invertible submatrix of $A$. In particular, for any $x$ satisfying $\|A^T(Ax-y)\|_\infty\le \sigma$, there is a solution $x^*$ of the \ref{pr: Dantzig Selector} approximating $x$ with error
\[
\|x-x^*\|_2\le \gamma\{2\sigma_k(x)_1+c_A\sigma+c_A\|A^T(Ax-y)\|_\infty\}
\le 
2\gamma \{\sigma_k(x)_1+c_A\sigma\}.
\]
\end{theorem}

%

However the choice of $p=1$ is the most interesting case, as $\ell_1$-norm is the closest convex norm to $\ell_0$-quasi-norm \cite{Ramirez_2013} and convex optimization is extremely well-nourished,  the shape of a unit ball associated with $\ell_p$-norm for  $0<p<1$ motivates researchers to explore the following nonconvex problem:
\begin{equation*}   \label{pr:nonconvex_relaxation}
 \min_{  x \in \mathbb R^N} \
\| {x}\|^p_p
\quad \text{subject to} 
\quad    Ax=y.
\end{equation*}
For certain values of  $p$ where   $0<p<1$, the above scheme leads to  more robust and stable theoretical guarantees compared to the case of $p=1$ \cite{Chartrand_2009,Saab_2008, wen2015stable}. In other words, much less restrictive recoverability conditions are achievable for certain $p$'s with  $0<p<1$ compared to $\ell_1$-norm recovery \cite{zheng2017does}. For instance, the following result implies that a sufficient condition for recovering an $s$-sparse vector in the noiseless case via $\ell_{0.5}$ minimization is $\delta_{3s}+27 \delta_{4s} < 26$, where an analogous result for $\ell_1$-norm recovery requires $\delta_{2s}+2\delta_{3s}<1$. We bring the following related theorem, which discusses compressible vectors. Since it is practically stringent to impose sparsity, it is desired to study compressible vectors.
\begin{theorem} \cite{Saab_2008}
Assume that $A\in \mathbb R^{m \times N}$ satisfies 
\[
\delta_{ks}(A)+k^{\frac{2}{p}-1}\delta_{(k+1)s}(A)< k^{\frac{2}{p}-1}-1,
\]
for $k>1$ such that $ks$ is a natural number. Given  $x\in \mathbb R^N$, let $y=Ax+\textbf{e}$ with $\|\textbf{e}\|_2\le \epsilon$. Then, an optimal solution $x^*$ of 
\begin{equation}  \label{l_p-minimization}
\min_{x\in \mathbb R^N} \|x\|_p^p \quad \mbox{ subject to } \quad \|Ax-y\|_2\le \epsilon,
\end{equation}
obeys  
\[
\|x^*-x\|_2^p\le C\epsilon^p+D s^{\frac{p}{2}-1} [{\sigma_s(x)}_p]^p,
\]
where 
\[
C  =2^p\Bigg[ \frac{1+k^{p/2-1}(2/p-1)^{-p/2}}
{ (1-\delta_{(k+1)s})^{\frac{p}{2}}   - (1+\delta_{ks})^{p/2}k^{p/2-1}} \Bigg],
\ \
D  =\frac{2(\frac{p}{2-p})^{p/2}}{k^{1-p/2}}
\Bigg[ 1+
\frac{  (1+k^{p/2-1})(1+\delta_{ks})^{p/2} }{ (1-\delta_{(k+1)s})^{p/2}- 
(1+\delta_{ks})^{p/2}k^{p/2-1}}
\Bigg].
\]
In particular, if $x$ is an $s$-sparse vector, then $\|x^*-x\|^p_2\le C\epsilon^p$.
\end{theorem}

The main issue with this result is its demand for a global minimizer of a nonconvex function for which there is no known theoretical guarantee. One way to bypass this is by utilizing classical schemes that obtain a local minimizer  with an initial point that is sufficiently close to a global minimizer. For example, the solution of least-squares is experimentally suggested  for this initialization \cite{Chartrand_2007}, while  there is no guarantee that such solution is in general close enough to a global minimizer. 
Since solving an unconstrained problem is simpler, another approach for taking advantage of nonconvex $\ell_p$-minimization (\ref{l_p-minimization}) with $0<p<1$  is tackling the following $\ell_p$-regularized least squares  problem:
\begin{equation}  \label{l_p-regularized}
\min_{x \in \mathbb R ^N} \ \frac{1}{2} \|Ax-y\|^2_2+\lambda \|x\|_p^p,
\end{equation}
where $\lambda>0$ is the  penalty parameter. Because  problems (\ref{l_p-minimization}) and (\ref{l_p-regularized}) are equivalent in limit, this parameter must be selected meticulously to obtain an approximate  solution for the original problem. Once $\lambda>0$ is fixed, algorithms such as iteratively-reweighted least squares explained in \cite{Scales_1988} are beneficial. 
Furthermore, the choice of $p$ plays a key role in the efficiency of this approach; for example, the best choices in image deconvolution are $p=1/2$ and $p=2/3$. The case of $p=1/2$ is a critical choice because it provides sparser solutions among  $p\in [1/2,1)$ and any $p\in (0,1/2)$ does not show significantly better performance \cite{Xu_Huo_2012}. As a result, it is essential to have a specialized algorithm for this choice \cite{Xu_2012}.

Despite  $p>1$ leads to strictly convex programming, it is theoretically shown that in this case not only the solution is not sparse, but also each entry is almost always nonzero. The same result holds  for all the primary problems mentioned above. 
\begin{theorem}\cite{Shen_2018}
Let $p > 1, N \ge m, \lambda>0, \mbox{and } \tau>0$. For almost all $(A, y)\in \mathbb{R}^{m \times N}$, a unique optimal solution $x^*_{(A,y)}$ to the any of 
the problems \ref{pr: Basis_Pursuit}, \ref{pr: Quadratically_Constrained_Basis_Pursuit}, \ref{pr: LASSO} and \ref{pr: Basis Pursuit Denoising},
when $\|.\|_1$ is replaced by $\|.\|_p^p$, satisfies $|\supp(x^*_{(A,y)})|=N$.
\end{theorem}
Recent applications in image processing, statistics and data science lead  to another class of sparse optimization problems subject to the nonnegative constraint, i.e.,
\begin{equation} \tag{NCS}  \label{pr: norm0-nonnegative-cs-problem}
 \min_{  x \in \mathbb R^N} \
\| {x}\|_0
\quad \text{subject to} 
\quad    A {x}= {y}, \quad   x \ge   0.
\end{equation}
This problem is similar to \ref{pr: initial_cs_problem} but their optimums  differ in general.  Similar convex/nonconvex relaxation approaches are advantageous to tackle \ref{pr: norm0-nonnegative-cs-problem} problem. Its convex relaxation leads to the so called Nonnegative Basis Pursuit (NBP) \cite{DonohoTanner_PNAS2005, KhajehnejadHassibi_TSP2011, Zhao_JORSC14}:
\begin{equation} \tag{NBP}  \label{pr: norm1-nonnegative-cs-problem}
\min_{  x \in \mathbb R^N} \
\| {x}\|_1
\quad \text{subject to} 
\quad    A {x}= {y}, \quad   x \ge   0.
\end{equation}
\begin{theorem} \cite{Zhao_JORSC14}
Any $x \ge 0$ such that $\|x\|_0\le s$ is recovered by \ref{pr: norm1-nonnegative-cs-problem} with $y=Ax$ if and only if for any index subset $S\subseteq [N]$ such that $|S|\le s$ there exists $\eta \in \mathcal{R}(A^T)$ such that $\eta_i=1$ {for $i \in S$} and $\|\eta_{S^c}\|_\infty<1$.
\end{theorem}
For a theoretical study on the following nonnegative version of LASSO  formulation, see \cite{itoh2016perfect}:
\begin{equation*}  \tag{NLASSO} \label{pr: NLASSO}
 \min_{  x \in \mathbb R^N }\
\quad  \frac{1}{2}\|  y-   A {x}\|^2_2+\lambda \| {x}\|_1 
\quad \text{subject to} 
\quad  x\ge 0.
\end{equation*}
Nonnegativity is an example of possible available prior information in the recovery process, which enables us to employ sparse measurement  matrices for the recovery of remarkably larger signals \cite{KhajehnejadHassibi_TSP2011}. 
In applications like infrared absorption spectroscopy, the non-zero elements of the original sparse signal are bounded. Imposing this boundary condition reduces the optimal set, and therefore its convex relaxation leads to a bounded/boxed basis pursuit that obtains better recovery properties \cite{donoho2010counting, donoho2010precise, liu2018binary}. 

Majority of studies in compressive sensing focuses on continues signals where $x \in \mathbb{R}^N$, even though there are several real-world applications that include discrete signals \cite{liu2018binary}. Examples are discrete control signal design, black-and-white or gray-scale sparse image reconstruction, machine-type multi-user communications, and blind estimation in digital communications. Here, we aim to reconstruct a discrete signal whose elements take their values from a finite set of alphabet. It is worth highlighting that when $x$ is discrete, even the $\ell_1$-norm recovery method leads to an NP-hard problem \cite{lange2016sparse}. Nevertheless, the study \cite{lange2016sparse} reveals this additional prior knowledge about the original signal  enhances the performance of both $\ell_0$ and $\ell_1$ norm recovery methods  via imposing new constraints.
For example, the sum of absolute values is designed to recover discrete sparse signals whose non-zero elements are generated from a finite set of alphabets with a known probability distribution. In binary compressive sensing, we can  represent this model as follows:
\begin{equation}  \label{SAV}
 \min_{  x \in \mathbb R^N }\
\quad  (1-p)\| x \|_1 + p \, \| x -  \textbf{1}_N \|_1
\quad \text{subject to} 
\quad  Ax=y,
\end{equation}
where $p$ is the probability of each component being one; analogous to the sparsity rate 
\cite{keiper2017compressed, liu2018binary,  nagahara2015discrete}. In this approach, when $p$ goes to zero ($x$ is sparse), the problem (\ref{SAV}) behaves like the \ref{pr: Basis_Pursuit}. The sum of norms represents the problem of binary compressive sensing as follows: 
\begin{equation*}  
 \min_{  x \in \mathbb R^N }\
\quad \|x\|_1 + \lambda \, \| x - \tfrac{1}{2} \, \textbf{1}_N \|_{\infty} 
\quad \text{subject to} 
\quad  Ax=y,
\end{equation*}
where the parameter $\lambda>0$ keeps $\ell_1$ and $\ell_{\infty}$ balanced \cite{wang2013binary}.
It is observed that the $\ell_{\infty}$-norm minimization 
tends to capture a representation whose coefficients have roughly the same absolute value \cite{fuchs,studer}.
Hence, the objective function
 $\|x\|_1 + \lambda \, \| x - \tfrac{1}{2} \, \textbf{1}_N \|_{\infty} $ is two fold. First, the sparsity is imposed  via  the $\ell_1$-norm term. Second, for those remaining   coefficients that are  deviated from 0, the binary  property  is  explored through the $\ell_{\infty}$ term. By the observation that solution of  the $\ell_{\infty}$-norm minimization favors to achieve  a solution with the  components of the same magnitude,  the $\ell_{\infty}$ term encourages those (deviated from 0) nonzero coefficients  to be centered at 1 (because $|x_i-\frac{1}{2}|$ obtains the same value for $x_i=0$ or 1). 
Nonetheless, the smoothed $\ell_0$ gradient descent technique has demonstrated to outperform such approaches in terms of recovery rate and time \cite{liu2018binary}.

 Finally, we emphasize on two important points. First, that most  results in sparse recovery are based on  asymptotic properties of random matrices but sometimes an application specifies a deterministic measurement matrix. Thus,  it is vital to have  tractable schemes  for testing  properties like the NSP, RIP, and etc. in  potential situations; see, e.g., \cite{d2011testing,tillmann2018computing}. In fact, the nature of this requirement,  encourages RIPless conditions \cite{candes2011probabilistic, zhang2013theory}. Second, all the introduced convex minimization problems can be solved by general purpose interior-point methods, however specified algorithms designed for these problems exist as well; see, e.g., \cite{burger2013adaptive,chambolle2011first,salahi, zhang2011unified}.

\section{Greedy Algorithms: Main Variations and Results}   \label{sec:greedy}

Greedy algorithms are  iterative approaches that take local optimal decisions in each step to eventually obtain a global solution. Hence, their success is due to some conditions on problem parameters. Greedy algorithms are mostly simple and fast, which  find numerous applications in various contemporary fields, including biology, applied mathematics, engineering, and  operations research. This supports the emerging interest in their performance analysis. 

Greedy algorithms   are efficient  in tackling those problems in compressive sensing as well.  Finding an $s$-sparse solution for an underdetermined linear system casts as the following  sparsity constrained problem:
\begin{equation} \tag{SC}  \label{pr: sparsity-constrained-problem}
\min_{  x \in \mathbb{R}^N} \
\| Ax-y\|_2^2
\quad \text{subject to} 
\quad    \|x\|_0\le s.
\end{equation}
The most popular greedy algorithm to solve this problem is the Orthogonal Matching Pursuit (OMP).
To tackle the problem \ref{pr: sparsity-constrained-problem}, starting from zero as the initial iteration, the following OMP  algorithm picks an appropriate index in each step to add its current support set and it  estimates the new iteration as the orthogonal projection of the measurement vector onto the subspace generated by corresponding columns of the current  support set \cite{Tropp_TIT2007,wang2015upport}. 
\begin{table}[h!]  \label{algorithm: OMP}
\noindent
\makebox[\textwidth]{
\begin{tabular}{l l }
\multicolumn{2}{c}{\textbf{Orthogonal Matching Pursuit}}\\  \hline \hline \\ [-2.6ex] 
Input: & $A\in \mathbb R^{m \times N}, y\in \mathbb R^m$, and initialize with $x_0=0\in \mathbb R^N, \mbox{ and } S_0=\emptyset$. \\ 
Iteration:& repeat until convergence:\\
 & $S_{n+1}$ = $S_n\cup \{j_{n+1}\} \quad \mbox{ such that } \quad j_{n+1}= \argmax_{j\in {S^c_n}} \big|\langle y-Ax_n,A_j\rangle\big|$,\\
 & $  x_{n+1} = \argmin \|  A   x-y\|^2_2 \quad \mbox{ subject to }   \quad \supp(x) \subseteq S_{n+1}$. \\
\hline \hline
\end{tabular}}
\end{table}

There is an extensive literature on the capability of the OMP in identifying the exact support set either in at most $s$ steps or at arbitrary many steps \cite{Tropp_TIT2004,Tropp_TIT2007, Cai_TIT2011}. 
\begin{theorem} \cite{Mo_2012, Wen_2013} \label{omp_recovery}
Assume that $A\in \mathbb{R}^{m \times N}$ satisfies the RIP of order $s+1$ such that $\delta_{s+1}\le \frac{1}{\sqrt{s}+1}$. Then, the OMP recovers any $s$-sparse vector $x\in \mathbb R ^N$ using $y=Ax$ in at most $s$ iterations.
\end{theorem}
Further, it is  proved that, when $s\ge 2$, for any $\delta$ such that $ 1/\sqrt{s+1}\le \delta<1$, any matrix with $\delta_{s+1}=\delta$ fails to recover all the $s$-sparse vectors \cite{dan2013sharp, Mo_2015}. It is still an open question that in terms of the OMP performance what happens if  $\delta_{s+1}\in (\frac{1}{\sqrt{s}+1}, \frac{ 1}{\sqrt{s+1}})$.
A result with the same spirit of Theorem \ref{omp_recovery} for noisy measurements requires an extra assumption on the magnitude of the smallest nonzero entry of the desired vector \cite{Wen_2017}. Here, we bring the next theorem that assumes  $\delta_s+ 2\delta_{31s}\le 1$. This inequality holds under the condition $\delta_{31s}\le 1/3$, which implies  the number of required measurements to successfully recover any $s$-sparse signal  via the OMP is $\mathcal{O}(s\ln N)$, that is, the number of measurements $m$ must be in the worst case  linear in sparsity level $s$. 
\begin{theorem} \cite{Zhang_2011}
Assume that  $A\in \mathbb{R}^{m\times N}$ satisfies 
$\delta_s+ 2\delta_{31s}\le 1$. Then, given an $s$-sparse vector $x\in \mathbb R^N$, the OMP obeys 
\[
\big \|x^{(30s)}-x\big \|_2\le \frac{2\sqrt{6}(1+\delta_{31s})^{\frac{1}{2}}}{1-\delta_{31s}}\|Ax-y\|_2.
\]
\end{theorem}
There are also results for the success of this algorithm in case of availability of partial information on the optimal support set   \cite{ge2019optimal, herzet2013exact}. This partial information  is  of the form of  either a subset of
the optimal support or an approximate subset with possibly wrong indices.

In many practical situations, there is no prior information available about the sparsity of a desired signal but generalized versions of the OMP  are still enabled to recover it  under a more stringent condition \cite{do2008sparsity}.

The OMP is also beneficial in recovering block sparse signals  \cite{eldar2010block, huang2019sharp,wen2018optimal}. To elaborate on this block version of the OMP, assume that $x=[x[1]^T,x[2]^T,\dots, x[L]^T]^T$ with $x[i]=[x_{d(i-1)+1},  x_{d(i-1)+2}, \dots, x_{di}]^T,$ $1\le i\le L$ (we assume the number of nonzero elements in different blocks are equal only for the sake of simplicity of exposition). Then, $x$ is called a block $s$-sparse signal if there are at most $s$ blocks (indices) $i$ such that $x[i] =0 \in \mathbb{R}^d$ \cite{eldar2010block}.  The support of this block sparse vector is defined as $\Omega:=\{i\, | \, x[i]\ne 0\}$. This definition reduces to the standard sparsity definition for $d=1$. 
We also partition  the measurement matrix $A$ as $A = \big[A[1],A[2], \dots ,A[L]\big]$, where  $A[i] =\big[A_{d(i-1)+1}, A_{d(i-1)+2}, \dots, A_{di}\big], \ 1\le i\le L$. Then, for a set of indices  $S=\{i_1,i_2,\dots, i_{|S|}\}$, we have $x[S]:= [x[i_1]^T,x[i_2]^T,\dots, x[i_{|S|}]^T]^T$, and $A[S] := \big[A[i_1],A[i_2], \dots ,A[i_{|S|}]\big]$. The block orthogonal matching  pursuit (BOMP) is as follows.

\begin{table}[h!]  \label{algorithm: BOMP}
\noindent
\makebox[\textwidth]{
\begin{tabular}{l l }
\multicolumn{2}{c}{\textbf{Block Orthogonal Matching Pursuit }}\\  \hline \hline \\ [-2.6ex]
Input: & $ A\in \mathbb R^{m \times Ld}, y\in \mathbb R^m$, and initialize with $x_0=0 \in \mathbb R^{Ld}, \mbox{ and } S_0=\emptyset$. \\ 
Iteration:& repeat until convergence:\\
 & $S_{n+1}$ = $S_n\cup \{j_{n+1}\}\quad  \mbox{ such that } \quad j_{n+1}= \argmax_{j\in [L]} \big\|A[j]^T (y-A[S_{n}]x[S_n] )\big \|_2$,\\
 & $  {x}[S_{n+1}] = \argmin_{x\in \mathbb{R}^{|S_{n+1}|d}}\big \|  A[S_{n+1}]   x-y \big \|^2_2$.\\ 
\hline \hline
\end{tabular}}
\end{table}
To present the exact recovery result of this algorithm, we need two other tools. 
\begin{definition} \cite{eldar2010block}
Consider a block measurement matrix $A$ introduced above. The block mutual coherence of $A$ is defined as 
\[
 \mu_B(A):=\max_{1\le l\ne r\le L} \ \frac{1}{d} \big \|A[l]^TA[r]\big \|_2,
\]
and its block subcoherence  is defined as 
\[
 \nu_B(A):=\max_{1\le l\le L}\ \max_{1\le i\ne j\le d} \ \frac{1}{d}\big |A_i[l]^TA_j[l]\big |,
\]
where $A_i[l]$ and  $A_j[l]$ are  the $i$th and $j$th columns of the  block $A[l]$.
\end{definition}
\begin{theorem} \cite{wen2018optimal}
Let $x$ be a block  sparse vector supported on $\Omega$  and  $y=Ax+\textbf{e}$ where $\|\textbf{e}\|_2\le \eta$ and the measurement matrix  satisfies $(2s-1)d \mu_B(A)+(d-1)\nu_B(A)<1$. Further, suppose that 
\[
\min_{i\in \Omega} \big\|x[i]\big\|_2  \ge \, \frac{   \eta \sqrt{2(1+(2d-1)\nu_B(A))} }{1-(2|\Omega|-1)d\mu_B(A)-(d-1)\nu_B(A)}.
\]
Then, the BOMP algorithm recovers block sparse vector $x$ in $|\Omega| $ iterations.
\end{theorem}
A similar result for the capability of another extended version of the OMP in recovering joint block sparse matrices  is found in \cite{fu2014block, shi2019sparse}. Despite many interesting features of the OMP, its  index selection is problematic. Precisely, if a wrong index is chosen, the OMP fails to expel this index so that the exact support cannot be found in this situation. The  Compressive Sampling Matching Pursuit (CoSaPM) is an algorithm to resolve this drawback \cite{needell2008iterative}. To find an $s$-sparse feasible vector for an underdetermined linear system, this algorithm allows $2s$ best potential indices enter the current support set. Then, the following CoSaMP  keeps $s$ entries that play the key role in the pertaining projection in the sense that their corresponding entries have most magnitude. 

\begin{table}[h!]  \label{algorithm: CoSaMP}
\noindent
\makebox[\textwidth]{
\begin{tabular}{l l }
\multicolumn{2}{c}{\textbf{Compressive Sampling Matching Pursuit }}\\  \hline \hline \\ [-2.6ex]
Input:&  $A\in \mathbb R^{m\times N}, y\in \mathbb R^m$, and initialize with $x_0=0\in \mathbb R^N, \mbox{ and } S_0=\emptyset$.\\
Iteration:&repeat until convergence:\\
 & $U_{n+1}$ = $\supp(x_n)\cup L_{2s}(A^T( y-Ax_n))$,\\
 & $  u_{n+1} = \argmin \|  A   x-y\|^2_2 \quad \mbox{ subject to } \quad  \supp(x) \subseteq U_{n+1}$,\\
 & $x_{n+1}=H_s(u_{n+1})$. \\

\hline \hline
\end{tabular}}
\end{table}
In the noiseless measurement case, the following theorem sates that any $s$-sparse signal, is recovered as the limit point of a sequence generated by the CoSaMP. 
\begin{theorem} \cite{FoucartRauhut_book2013}
Assume that $A\in \mathbb{R}^{m \times N}$ satisfies $\delta_{8s}<0.4782$. Then, for any $x\in \mathbb{R}^{N}$ and $\textbf{e}\in \mathbb R^m$, the sequence $x_n$  generated by CoSaMP using $y=Ax+\textbf{e}$, where $s$ is replaced by $2s$,  obeys 
\[ 
 \|x-x_n\|_2\le Cs^{-1/2}{\sigma_s(x)}_1  +D \|\textbf{e}\|_2+2 \rho ^n\|x\|_2,
\]
where constants $C,D$ and $\rho\in (0,1)$ only depend on $\delta_{8s}$. In particular, if $\tilde{x}$ denotes a cluster point of this sequence, then
\[ 
 \|x-\tilde{x}\|_2 \le  Cs^{-1/2}{\sigma_s(x)}_1 +D \|\textbf{e}\|_2.
\]
\end{theorem}
The following generalized orthogonal matching pursuit (gOMP) allows a given number of indices $t\ge 1$ enter a current support set. Consequently, a faster exact  recovery under the RI constants involved conditions  is attained \cite{WangKS_TSP2012,wen2017novel}.
\begin{table}[h!]  \label{algorithm: gOMP}
\noindent
\makebox[\textwidth]{
\begin{tabular}{l l }
\multicolumn{2}{c}{\textbf{Generalized Orthogonal Matching Pursuit }}\\  \hline \hline \\ [-2.6ex]
Input:&  $ A\in \mathbb R^{m\times N}, y\in \mathbb R^m, t\in \mathbb N$, and initialize with $x_0=0\in \mathbb R^N, \mbox{ and } S_0=\emptyset$.\\
Iteration:&repeat until convergence:\\
 & $U_{n+1}$ = $U_n\cup L_{t}\big(A^T( y-Ax_n)\big)$,\\
 & $  x_{n+1} = \argmin \|  A   x-y\|^2_2 \quad \mbox{ subject to } \quad  \supp(x) \subseteq U_{n+1}$.\\

\hline \hline
\end{tabular}}
\end{table}
A recent algorithm combines the BOMP and the gOMP to propose a block generalized orthogonal matching pursuit for recovering block sparse vectors with possibly different number of nonzero elements in each block \cite{yang2018recovery}.
A closely related task to the problem \ref{pr: sparsity-constrained-problem} is the following  nonnegative sparsity constrained (NSC) problem 
\begin{equation} \tag{NSC}  \label{pr: nonnegative-sparsity-constrained-problem}
\min_{  x \in \mathbb{R^N}} \
\| Ax-y\|_2^2
\quad \text{subject to} 
\quad    \|x\|_0\le s \quad \mbox{and } \quad x\ge 0.
\end{equation}
Then, its solution may  completely differ from the one to \ref{pr: sparsity-constrained-problem}.
Bruckstein et al. \cite{Bruckstein&Elad&Zibulevsky_TIT2008}  presented an adapted version of the OMP  for finding nonnegative sparse vectors of an underdetermined system, namely, the Nonnegative Orthogonal Matching Pursuit (NOMP) below. They demonstrated its capability to find sufficiently sparse vectors.
In practice, it is probable not to have sparsity level in hand, this challenge is also doable \cite{Needell_2010}.

\begin{table}[h]  \label{algorithm: NOMP}
\noindent
\makebox[\textwidth]{
\begin{tabular}{l l }
\multicolumn{2}{c}{\textbf{Nonnegative Orthogonal Matching Pursuit}}\\  \hline \hline \\ [-2.6ex]
Input: & $A\in \mathbb R^{m \times N}, y\in \mathbb R^m$, and initialize with $x_0=0\in \mathbb R^N, \mbox{ and } S_0=\emptyset$. \\ 
Iteration:& repeat until convergence:\\
 & $S_{n+1}$ = $S_n\cup \{j_{n+1}\} \quad \mbox{ such that } \quad j_{n+1}= \argmax_{j\in {S^c_n}} \langle y-Ax_n,A_j\rangle_+$,\\
 & $  x_{n+1} = \argmin \|Ax-y\|^2_2 \quad \mbox{ subject to } \quad \supp(x) \subseteq S_{n+1} \mbox{ and } x\ge 0$. \\

\hline \hline
\end{tabular}}
\end{table}

There is an emerging interest to explore novel greedy algorithms to find a sparsest feasible point of a set \cite{Bahmani_2013, Beck_2015}. This casts as a minimization problem under sparsity constrained possibly intersecting a desired set. Nevertheless, the main idea of a greedy algorithm here is to start with a feasible sparse point (possibly zero) and add one or  several candidate indices to the current support set. Once the support set is updated, a new iteration is obtained via a projection step. The  Constrained Matching Pursuit (CMP) \cite{shen_2019} investigates the following constrained sparse problem:
\begin{equation*} \label{pr: constrained_sparse_problem}
\min_{  x \in \mathbb{R^N}} \
\| x\|_0
\quad \text{subject to} 
\quad    Ax=y \quad \mbox{ and } \quad x\in \mathcal P,
\end{equation*}
where $\mathcal P \subseteq \mathbb R ^N$ is a closed constraint set containing the origin. 
One can see that the OMP and NOMP can be employed  for solving the above problem with the constraint set $\mathcal P$ with $\mathbb R^N$ and $\mathbb R^N_+$, respectively.

\begin{table}[h]  \label{algorithm: CMP}
\noindent
\makebox[\textwidth]{
\begin{tabular}{l l }
\multicolumn{2}{c}{\textbf{Constrained Matching Pursuit }}\\  \hline \hline \\ [-2.6ex]
Input: & $A\in \mathbb R^{m \times N}, y\in \mathbb R^m, \mathcal P \subseteq \mathbb R^N$, and initialize with $x_0=0, \mbox{ and } S_0=\emptyset$. \\ 
Iteration:& repeat until convergence:\\
&  $S_{n+1}$ = $S_n\cup \{j_{n+1}\} \quad \mbox{ such that } \quad j_{n+1}= \argmin_{j\in {S^c_n}}  g^*_j $  \mbox{ with }\\
& $g^*_j=  \min_{t \in \mathbb R} \| y - A (x_n + t \, \mathbf e_j  ) \|^2_2 \quad \mbox{ subject to }  \quad x_n + t \, \mathbf e_j  \in \mathcal P$,\\

& $ x_{n+1} = \argmin \|Ax-y\|^2_2 \quad \mbox{ subject to } \quad \supp(x) \subseteq S_{n+1} \mbox{ and } x\in \mathcal P$. \\
 
\hline \hline
\end{tabular}}
\end{table}

The performance of the CMP relies on the measurement matrix $A$ as well as the involved constraint set $\mathcal P$. In fact, there exists a set $\mathcal P$ for which there is no measurement matrix $A$ such that the CMP is successful \cite{shen_2019}. Hence, a class of convex coordinate-projection (CP) admissible sets is of interest. A nonempty set $\mathcal P \subseteq \mathbb R^N$ is called CP admissible if for any $x \in \mathcal P$ and any index set $\Jcal \subseteq \supp(x)$, we have $x_{\Jcal}\in \mathcal P$. The conic hull structure of the CP sets can be used to extend the RIP and the RO constants over these sets. A verifiable   exact recovery condition is then developed based on such constants for the closed convex CP admissible cones \cite{shen_2019}. (Consequently, a verifiable condition for the NOMP is in hand because the nonnegative orthant is a closed convex CP admissible cone). These verifiable conditions are often NP-hard to check, like almost all the introduced tools in this survey such as the NSP, and  RIP. But random measurement matrices  with the right size  hold these verifiable conditions with a high probability. This fact  confirms the significance of random measurements
in sparse optimization  \cite{tillmann2013computational}.

There are many other effective algorithms in sparse optimization that are applicable  for the compressive sensing setting as well. For example, threshholding based algorithms use the adjacent matrix $A^*$ for approximating inversion action and exploit the hard thresholding operator to solve the square system $A^*Ax=A^*y$ via a fixed-point method \cite{blumensath2011iterative, carrillo2011iterative, daubechies2004iterative,   figueiredo2003algorithm, fornasier2008iterative}. 
Iterative  reweighted algorithms often outperform algorithms presented in Section \ref{sec:lp_recovery}. The iterative reweighted $\ell_p$-algorithm is as follows: 
\begin{equation*} \label{eqn:reweighted_l1}
\quad x_{n+1}\in \underset{x \in \mathbb R^N} \argmin \ \|W_nx\|_p \quad \text{subject to} \quad \|Ax-y\|_2\le \epsilon,
\end{equation*}
where $ W_n$ is a diagonal weight matrix defined based on the current iteration $x_n$. For example, $W^k=\text{diag}(w_n^1, w_n^2, \dots, w^N_n)$ where $w^j_n=1/(|x_n^j|+\gamma)^t$ with $t$ and $\gamma >0$, which encourages small elements  to approach zero quickly \cite{abdi2013cardinality, candes2008enhancing,chartrand2008iteratively}.
The so called smoothed $\ell_0$-norm takes advantage of  smooth approximations of the  $\ell_0$ quasi-norm  to tackle the \ref{pr: initial_cs_problem} and \ref{pr: norm0-nonnegative-cs-problem} \cite{Mohammadi_SP2014, Mohimani_TSP2009}. 

In traditional compressive sensing approaches, we recover sparse signals from $m = \mathcal O{\left({{s}{\log{N/s}}}\right)}$ linear measurements. Model-based compressive sensing enables us to substantially reduce the number of measurements to $m = \mathcal O{\left({s}\right)}$ without loosing the robustness of the recovery process \cite{baraniuk2008model}. In fact, model based compressive sensing includes the structural dependencies between the values and locations of the signal coefficients. For example, modeling the problem of binary compressive sensing with a bi-partite graph and representing the recovery process as an edge recovery scheme recovers binary signals more accurately, compared to binary $\ell_1$-norm and tree-based binary compressive sensing methods \cite{nakarmi2012bcs}. 

On the other hand, quantum annealers are type of adiabatic quantum computers that tackle computationally intensive problems, which are intractable in the realm of classical computing. From a problem solving point of view, quantum annealers receive coefficients of an Ising Hamiltonian as input and return a solution that minimizes the given energy function in a fraction of a second \cite{ayanzadeh2019quantum-assisted,ayanzadeh2020reinforcement,nishimori2017exponential}. Recent studies have revealed that well-posed binary compressive sensing and binary compressive sensing with matrix uncertainty problems are tractable in the realm of quantum computing \cite{ayanzadeh2020leveraging,ayanzadeh2019sat, ayanzadeh2020ensemble,ayanzadeh2019quantum}.

\section{Numerical Performance Analysis } \label{sec:numerical_results}
In this section, we numerically examine the performance of the most prevalent sparse recovery methods introduced  in this survey in recovering unigram text representation from their linear embedding (measurement). To describe the sparse recovery task of interest, we start by an introductory explanation on  document (text) unigram representations  and  embeddings. Suppose that a vocabulary $\mathcal V$ is available, namely, a set of all the considered words in a context. Then, a unigram representation of a document is simply a vector in $\mathbb R^{|\mathcal V|}$, where its $i$th entry counts the occurrences of the $i$th word of  $\mathcal V$ in that document  \cite{pang2004sentimental,paskpv2013compressive}. Unigram representations are highly sparse as the size of a vocabulary is often too large. These representations are expensive to store and further may not be able to effectively capture the  semantic relations among the words of a vocabulary, so text embeddings are naturally favorable.  

We start by defining word embeddings and discuss text embeddings afterward. The goal in word embedding is to encode each word of a vocabulary into a vector representation in a much lower dimensional space $m\ll |\mathcal V|$ such that the semantic relations among words are preserved. Word embeddings  have lately gained much attentions in broad applications of the natural language processing such as  classification, question answering and part of speech tagging.  There are different linear and nonlinear effective approaches for word embedding, for example,
pretrained models Global Vectors For Word Representation (GloVe) \cite{Glo}, word2vec \cite{mikolov2013distributed} and 
Rademacher embeddings. The GloVe is a neural (network) based word embedding, where  the word2vec is mainly based on matrix decomposition. Assuming that word embeddings are available for all the  words of a vocabulary, a text embedding is simply a linear combination of them with the coefficients coming from its unigram representation. In other words, we have $A x_T^{\mathrm{unigram}}=y_T^{\mathrm{embedding}}$,  where $A\in \mathbb R^{m \times | \mathcal V|}$ is an embedding (or measurement matrix) generated from a specific methodology,
$x_T^{\mathrm{unigram}}\in \mathbb R^{|\mathcal V|}$ is the unigram representation of a text, and $y_T^{\mathrm{embedding}}\in \mathbb R ^m$ is its (unigram) embedding. 
The measurement matrices in our experiments  are generated via GloVe embedding method and 
the Radamacher distribution with the following 
probability mass function:
\begin{equation*}
f(k)=
\begin{cases}
1/2& \mathrm{if}\,k=-1, \\
1/2 & \mathrm{if}\,k=+1.\\
\end{cases}
\end{equation*}
These matrices are:
$m\times 17000$ for MR movie reviews \cite{pang2005seeing} and $m\times 20000$ for 
SUBJ subjectivity dataset \cite{pang2004sentimental}, respectively. Embedding sizes  for our experiments are $m=50,100,200,300$ and $1600$.
\begin{figure}[h!]
\centering
\begin{subfigure}[b]{0.235\textwidth}
\includegraphics[width=\textwidth]{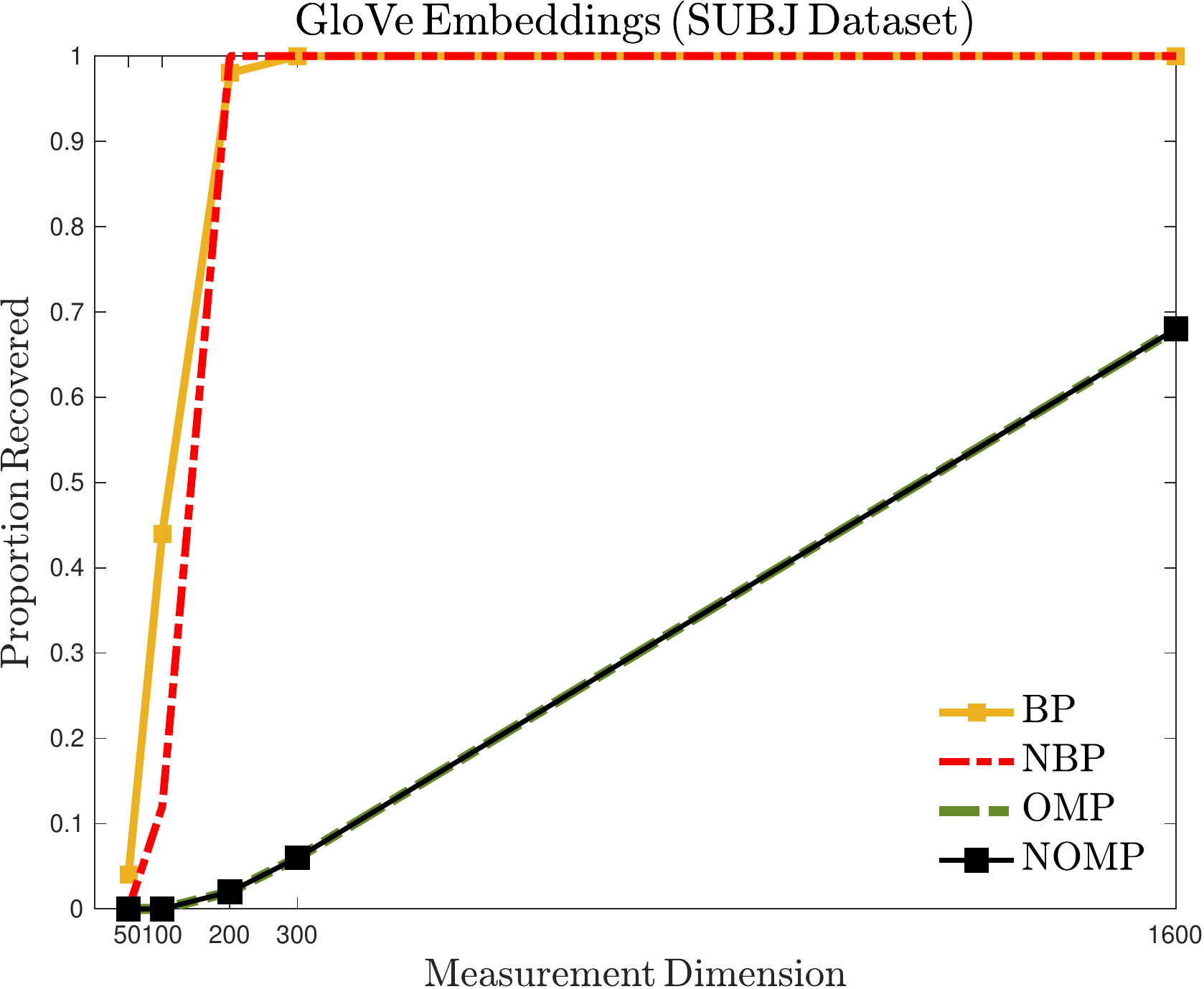}
\end{subfigure}
~ 
\begin{subfigure}[b]{0.235\textwidth}
\includegraphics[width=\textwidth]{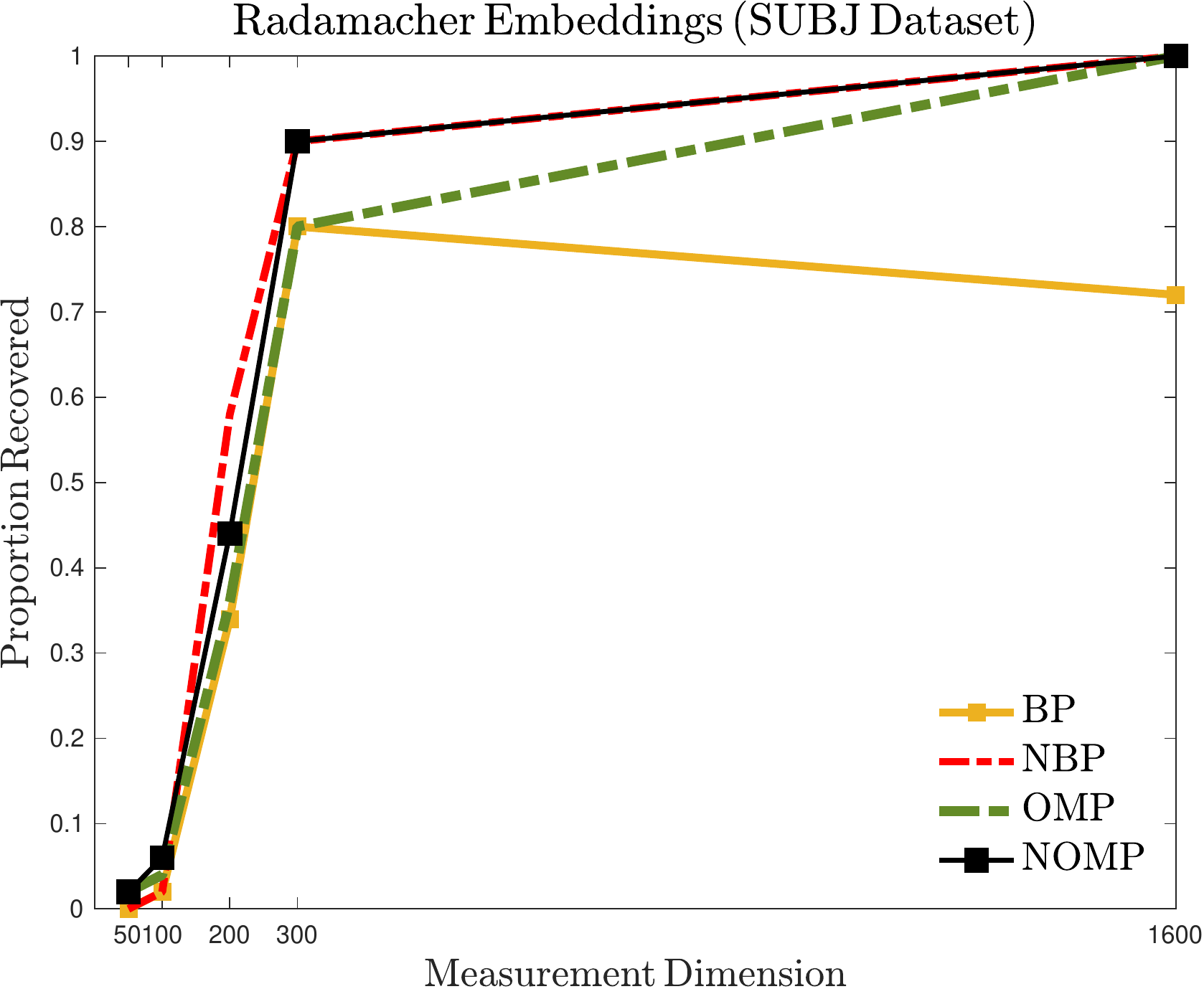}
\end{subfigure}
~ 
\begin{subfigure}[b]{0.235\textwidth}
\includegraphics[width=\textwidth]{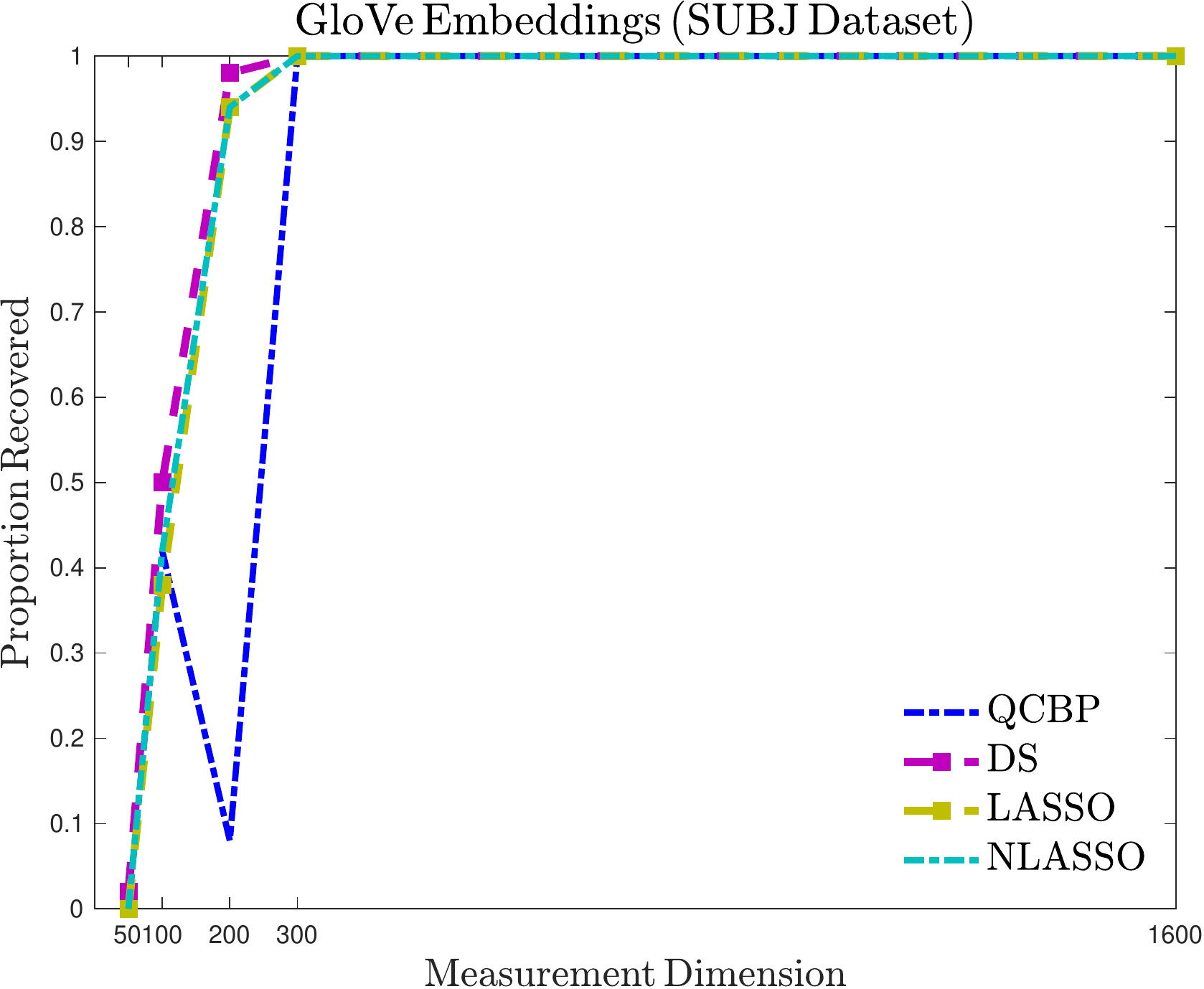}
\end{subfigure}
~ 
\begin{subfigure}[b]{0.235\textwidth}
\includegraphics[width=\textwidth]{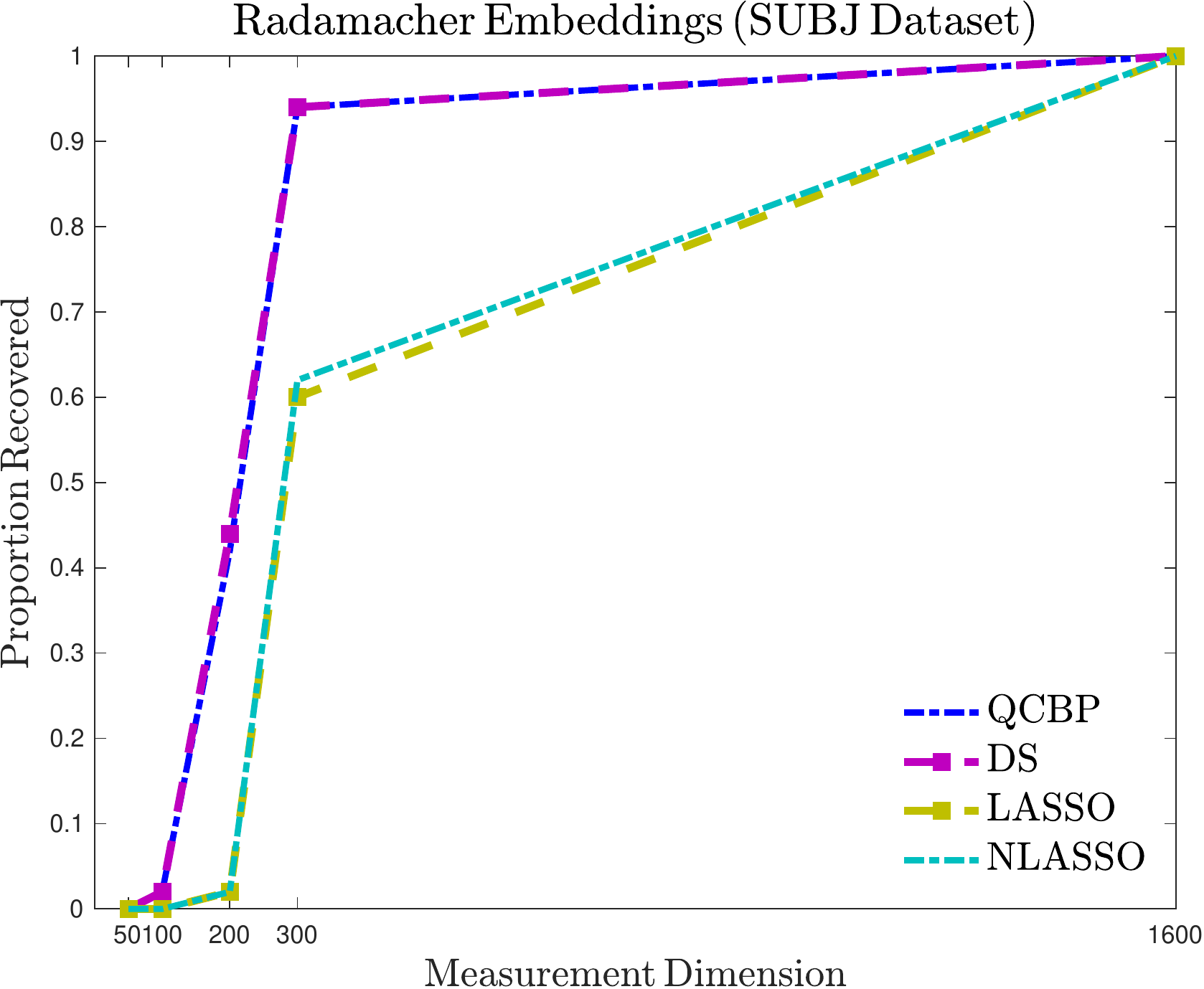}
\end{subfigure}  
\caption{SUBJ Dataset recovery}
\end{figure}

\begin{figure}[h!]
    \centering
    \begin{subfigure}[b]{0.235\textwidth}
    \includegraphics[width=\textwidth]{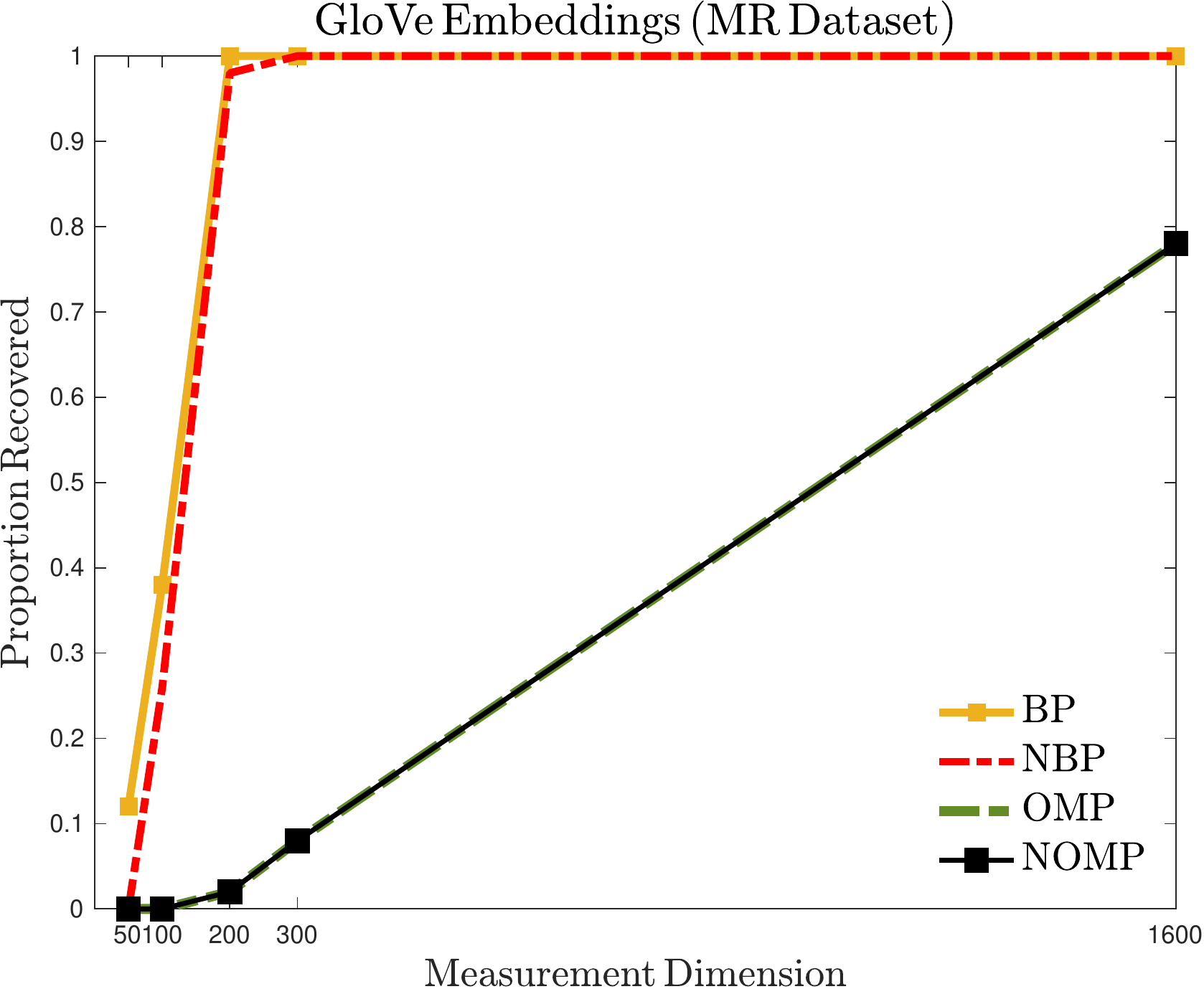}
    \end{subfigure}
    ~ 
    \begin{subfigure}[b]{0.235\textwidth}
    \includegraphics[width=\textwidth]{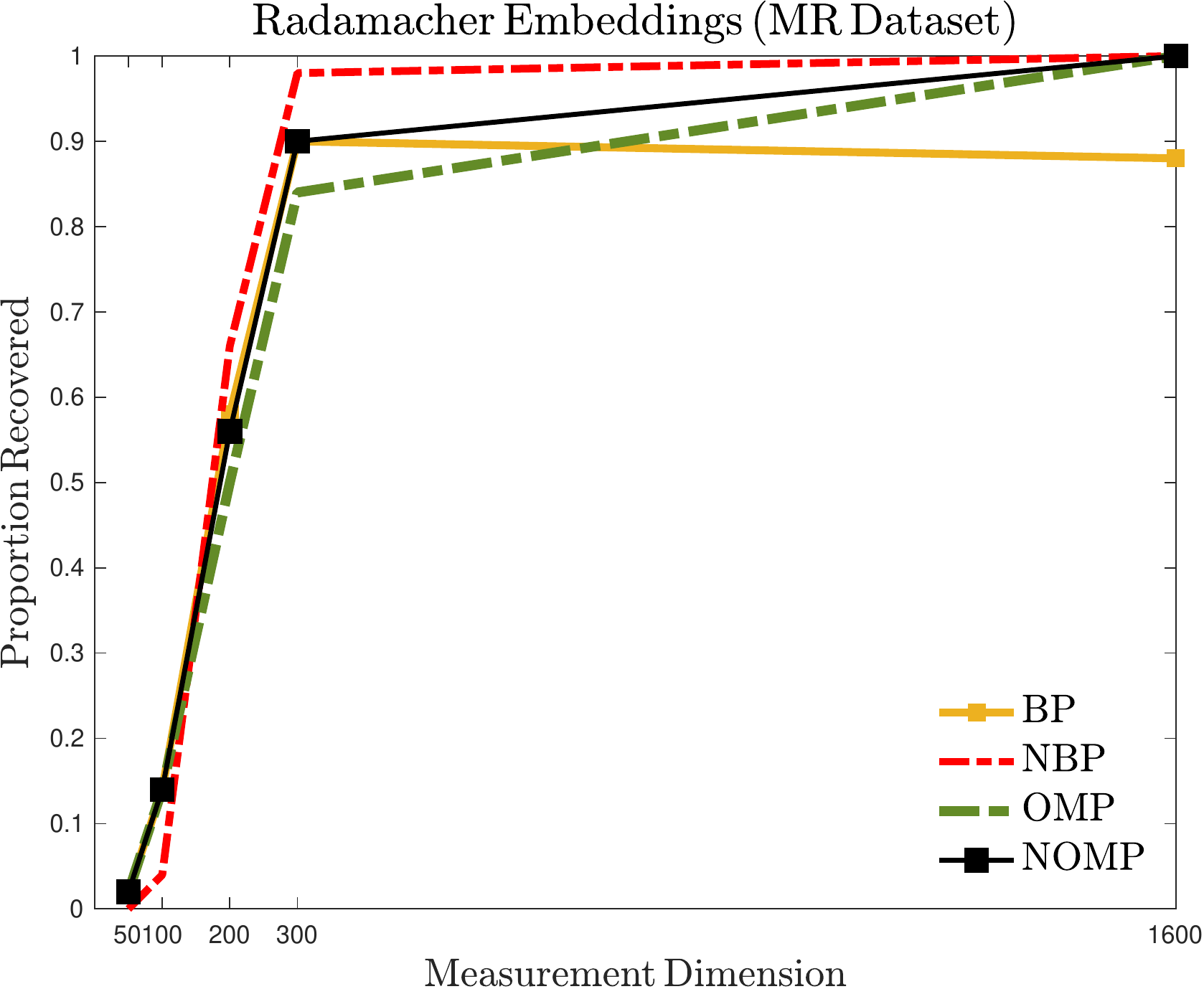}
    \end{subfigure}
    ~ 
    \begin{subfigure}[b]{0.235\textwidth}
    \includegraphics[width=\textwidth]{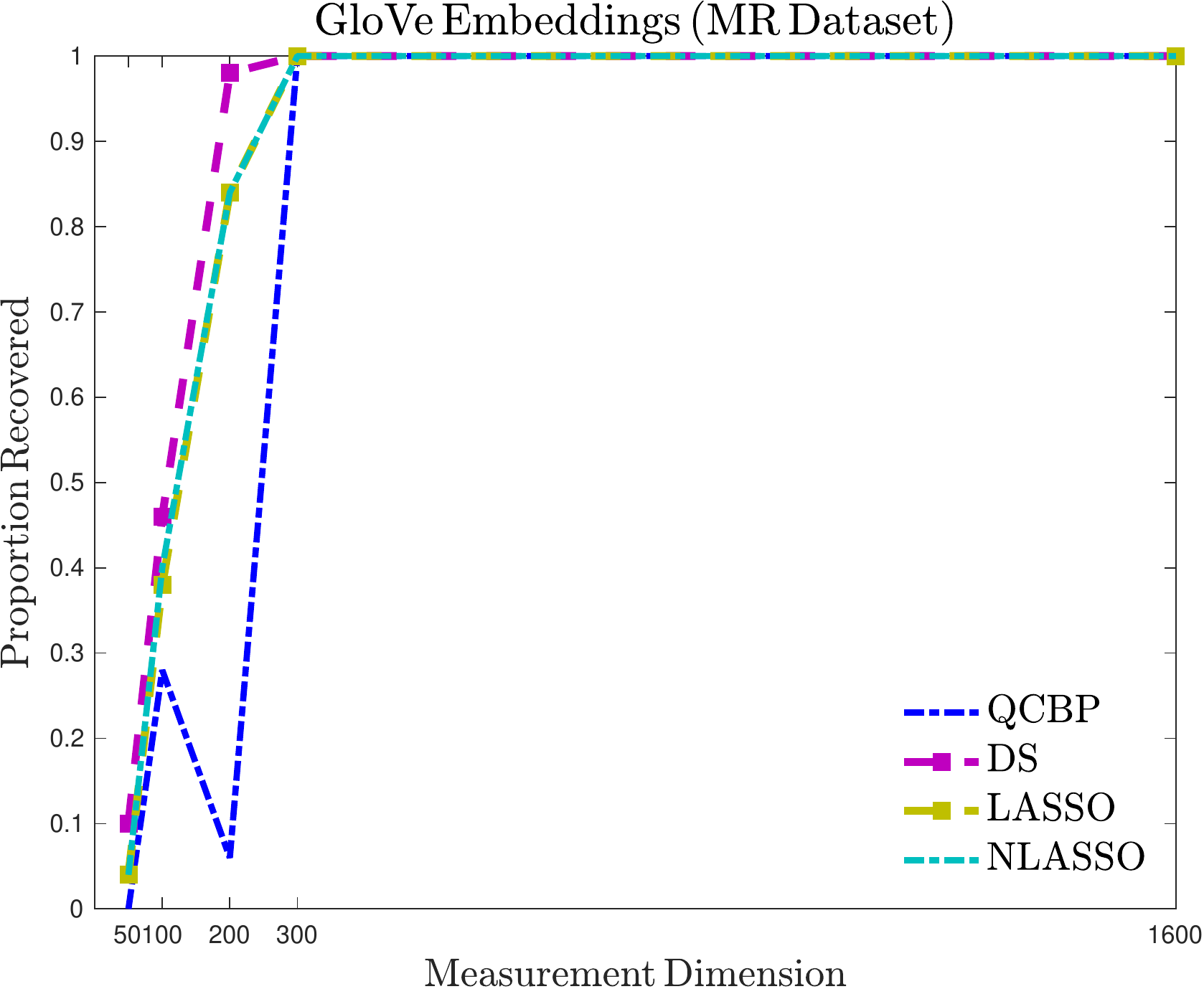}
    \end{subfigure}
    ~ 
    \begin{subfigure}[b]{0.235\textwidth}
    \includegraphics[width=\textwidth]{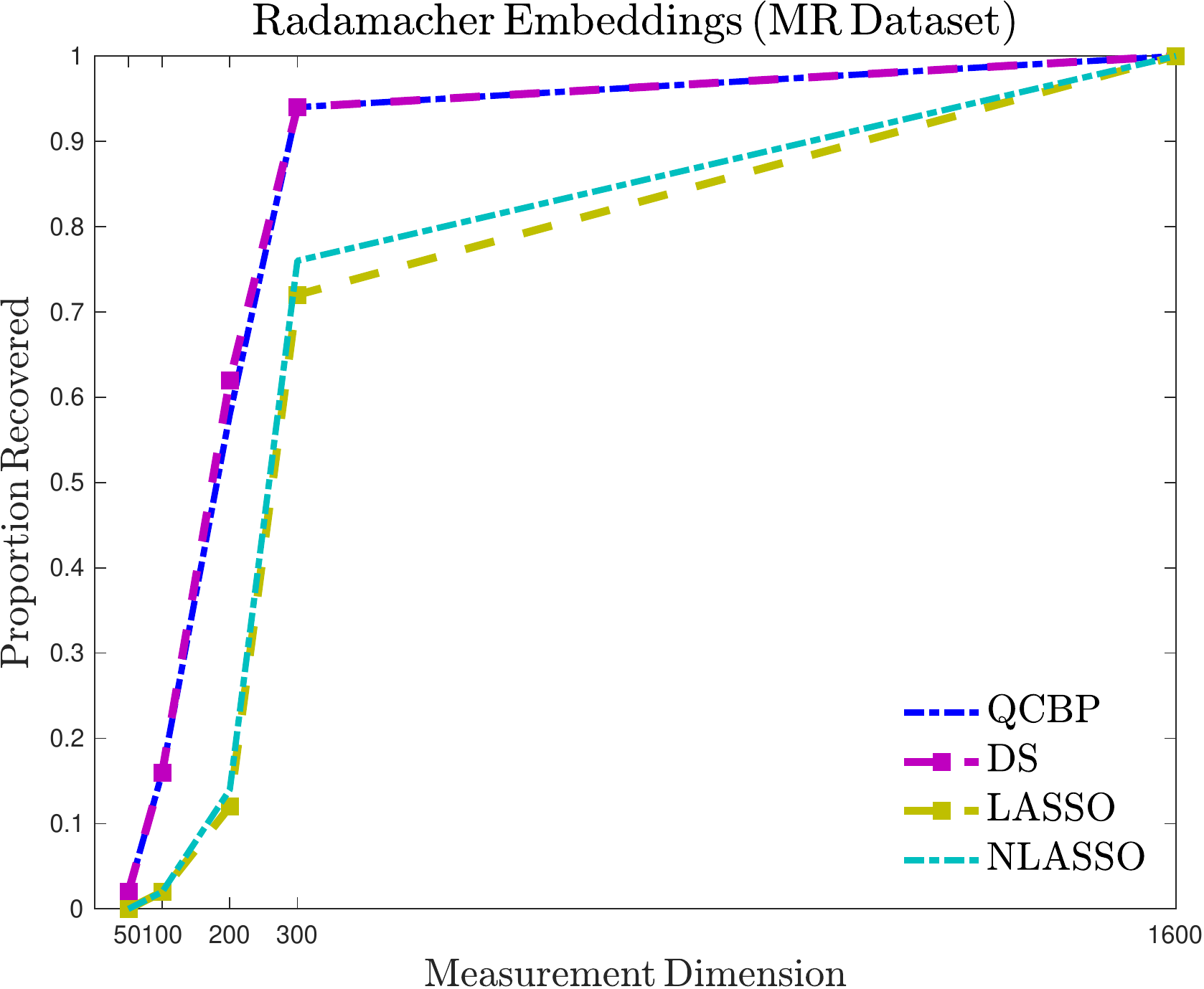}
    \end{subfigure}  
    \caption{MR Dataset recovery}
\end{figure} 

These figures confirm the efficiency of $\ell_1$ recovery and greedy algorithms in recovering unigram representations of $50$ documents, where success  is achieved if the relative error is smaller than $10^{-7}$.
We avoid a detailed explanation on the implemented algorithms  for this specific application rather finish this survey with a more general conclusion as follows. The $\ell_p$ recovery with $0<p\le 1$, where $p<1$ promises  better theoretical results than $p=1$, and greedy algorithms are both effective in recovering sparse vectors in various applications if the measurement matrix inherits several properties explained in this survey. In the constrained case, the constraint set also plays a crucial role in the sparse recovery conditions as well \cite{shen_2019}. Nevertheless, these properties often provide sufficient conditions and for some  applications specified measurement matrices can work properly as well and even better than random matrices \cite{arora2018compressed}.  Further, except when the sparsity level is relatively small, the $\ell_p$ recovery with $0<p\le 1$ is superior to the greedy approach. 
For a more comprehensive study on the computational complexities and efficiency of these algorithms or related studies, see, e.g., \cite{adhikari2017nonconvex, fornasier2010numerical, qiasar2013compressive,zhang2011user}.
\\

\textbf{Acknowledgement}
\\
The authors would like to greatly thank  Professors  Jinglai Shen and Maziar Salahi for their meticulous comments and constructive suggestions on this work.

\bibliographystyle{plain}
\bibliography{Survey}

\end{document}